\documentclass[12pt]{amsart}
\usepackage{amscd}
\advance\hoffset by -0.5 truecm
\usepackage{amssymb}

\def \0{\lambda_{0}}

\usepackage{graphicx, amssymb, amsfonts}

\usepackage{epsfig}

\begin{document}

\title{Affine Configurations and Pure Braids}
\author{Pablo Su\'arez Serrato}
\address{Girton College, Cambridge.}
\address{Mathematisches Institut LMU, Theresienstrasse 39, M\"unchen 80333 Deutschland.}
\email{Pablo.Suarez-Serrato@mathematik.uni-muenchen.de}
\thanks{This work was carried out at the Instituto de Matem\'aticas, UNAM and partially funded by the Sistema Nacional de Investigadores and a doctoral fellowship from CONACyT M\'exico. }

\begin{abstract} We show that the fundamental group of the space of ordered affine--equivalent configurations of at least five points in the real plane is isomorphic to the pure braid group modulo its centre. In the case of four points this fundamental group is free with eleven generators.
\end{abstract}

\maketitle

\section{Introduction}

The aim of this note is to establish a new link between discrete geometry and classical concepts in topology. Configurations of vectors can be understood with the help of Grassmannian manifolds equipped with a stratification controlled by matroids \cite{BS}. Similarly, affine configurations yield stratified Grassmannians whose structure is determined by affine oriented matroids. The decomposition provided by Schubert cells of a Grassmannian has striking equivalent formulations in terms of convex polyhedra and of matroids \cite{GGMS}. Affine analogues of these results have been pursued in \cite{ABM, BMO}, where the focus there lies in understanding the combinatorics of these polyhedral structures. We will show that the topology of based loops of certain affine configurations can be explained using braids.

A braid should be thought of as a system of $k$ strings which do not intersect and have fixed start and end points. Let $X_{k}$ and $Y_{k}$ be the sets $\{ (1, 0, 1), (2,0,1), \ldots, (k,0,1) \}$ and $\{ (1, 0,0), (2,0,0), \ldots, (k,0,0)\}$ in ${\bf R}^3$, respectively. Consider a smooth embedding $\beta_{i}$ of an interval $I$ into ${\bf R}^3$ such that $\beta_{i}(0)$ lies in $X_{k}$ and $\beta_{i}(1)$ in $Y_{k}$. Denote by $\beta$ a collection of $k$ such embeddings  $\beta_{i}$, each one starting from a distict point in $X_{k}$ and ending in a distinct point in $Y_{k}$, such that the projection to the $z$ coordinate is decreasing as $t$ runs from $0$ to $1$ in $I$. We also require that $\beta_{i}(t) \neq \beta_{j}(s)$  for every $t,s$ in $I$ and $i \neq j$, which means we obtain disjoint images and the strings do not intersect.

{\bf Definition} {\it A braid $b$ is an isotopy class of a $\beta$ as described above.}

Braids form a group  $B_{k}$ \cite{Ar} and each braid $b$ permutes the order of $X_{k}$ to a new order in $Y_{k}$. When this order remains unchanged, i.e. the permutation is trivial, we say $b$ is a {\it pure} braid. It is easy to see that pure braids form a subgroup $P_{k}$. The fundamental group of the space of $k$-tuples of distinct points in the plane is isomorphic to the pure braid group $P_k$ \cite{FN}. The centre of $P_k$ is the infinite cyclic group $Z_k$,  generated by the element $\Delta_{k}$.

The spaces of affine configurations of flats in affine and projective spaces were defined in \cite{ABM}. Following the same notation we consider $\overrightarrow{ \mathbb{A}^{2}_{k,0}}(\neq)$, the space of distinct $(\neq)$, oriented $(\rightarrow)$ affine classes of sets of $k$ points in the real affine plane. For simplicity we will refer to this space as $A_k$.

{\bf Theorem~1.} {\it  For $k\geq 5$, $\pi_1(A_k, \ast )\cong P_k / Z_k $.}

This can be seen as a variation of the classical result about the configurations of distinct points in the plane and the pure braid group mentioned above. It allows results from braid groups to be translated into results about affine configurations. For example, the conjugacy problem for braids is solvable \cite{Gar}, therefore the isomorphism above implies  that the conjugacy problem for homotopy classes of based loops of oriented affine configurations in the plane is also solvable. In principle there exists an algorithm to decide whether or not two such loops are conjugate.

Suppose $X$ is a topological space and let ${\rm Homeo}(X)$ denote the group of self-homemorphisms of $X$. Let  ${\rm Homeo}_0(X)$ and  ${\rm Homeo}^{+}(X)$ be the subgroups which consist of homeomorphisms which are isotopic to the identity and those which preserve orientation, respectively. Define the {\it pure mapping class group} of $X$ as  ${\rm Homeo}^{+}(X)/ {\rm Homeo}_0(X) $.

The group $P_k$ is isomorphic to the pure mapping class group of the disk $D^2$ minus $k$ points \cite{Bir}, we can now establish another isomorphism. From Theorem 4.5 in \cite[p.164]{Han} it follows that the pure mapping class group of the sphere $S^2$ with $k+1$ punctures is $P_k / Z_k$. Therefore we obtain for $k \geq 5$,

{\bf Corollary ~2.} {\it  The fundamental group of $A_k$ is isomorphic to the pure mapping class group of the sphere $S^2$ punctured $k+1$ times. }

We will show that the group  $\pi_1(A_4, \ast)$ is isomorphic to a free group with eleven generators. This can be seen by using  the combinatorial decomposition of the affine configuration spaces, as outlined in \cite{ABM} and \cite{BMO}. Notice that the configuration space $A_3$ consists of a single point.

{\bf Acknowledgements:} The author thanks Luis Montejano Peimbert for his guidance and encouragement, and is also grateful to Francisco Gonzalez Acu\~na, Mario Eudave Mu\~noz and Raymond Lickorish for useful conversations and for pointing out subtle details in previous versions.

\section{Proofs}
Let $ \langle p_1, p_2, ... , p_k\rangle $ be the affine span of the set $ \{ p_1, p_2, ... , p_k \} $, where $ p_i \in {\bf R}^n $. Define
$$ E_k:= \{ (p_1, p_2, ... , p_k) : p_i \in {\bf R}^2 ,\langle p_1, p_2, ... , p_k\rangle ={\bf R}^2  , p_i \neq p_j \} . $$

In other words, we are considering the sets of $k$ points in $\bf{R}^2$ which are not all co-linear and hence span the real plane affinely. The group ${ \rm{Aff} }_{+}({\bf R}^2)$ of orientation preserving affine transformations of ${\bf R}^2$ acts freely on $E_{k}$. We denote the quotient of $E_k$ under this action by $A_k:= E_k / { \rm{Aff} }_{+}({\bf R}^2)$, in \cite{ABM} it is denoted by $ \overrightarrow{ \mathbb{A}^{2}_{k,0}}(\neq) $. We call $A_k$ the space of oriented affine configurations of $n$ distinct points in the real plane. It was shown in \cite{BMO} that  ${\rm{Aff}}_{+}({\bf R}^2)\stackrel{i}{\rightarrow}  E_k \rightarrow A_k$ is a principal fibre bundle.
\par

{\bf Theorem~1.} {\it  For $k\geq 5$, $\pi_1(A_k, \ast )\cong P_k / Z_k $.}

\begin{proof} Consider the exact sequence of homotopy groups associated to the principal fibre bundle ${\rm{Aff}}_{+}({\bf R}^2)\stackrel{i}{\rightarrow}  E_k \rightarrow A_k$,
$$...\rightarrow \pi_1({\rm{Aff}}_{+}({\bf R}^2), \ast)\stackrel{i_{\ast}}\rightarrow \pi_1(E_k , \ast)\rightarrow \pi_1(A_k, \ast)\rightarrow \pi_0({\rm{Aff}}_{+}({\bf R}^2), \ast )\cong  1.$$

Note that $\pi_1({\rm {Aff}}_{+}({\bf R}^2), \ast)$ is isomorphic to ${\bf Z}$, because ${\rm {Aff}}_{+}({\bf R}^2)$ retracts to ${\rm SO}(2)$. Say the generator of $\pi_1({\rm {Aff}}_{+}({\bf R}^2), \ast)\cong {\bf Z}$ is $g$. Since group morphisms send generators to generators we have, $$\langle i_{\ast}(g) \rangle  =  i_{\ast}(\pi_1({\rm {Aff}}_{+}({\bf R}^2), \ast)).$$

Let $F_k := \{ (p_1, p_2, ... , p_k) : p_i \in {\bf R}^2 , p_i \neq p_j \}$ be the standard space of configurations of distinct ordered $k$-points in ${\bf R}^2$.

Notice that $g$, as a path in ${\rm {Aff}}_{+}({\bf R}^2)$, can be represented as a continuous full rotation of the plane.

Take $\gamma_{\theta}= \left( \begin{array}{cc} \cos \theta & \sin \theta \\  -\sin \theta  & \cos \theta \end{array} \right) $, for $\theta \in [0,2\pi]$ to be a representative of $g$. Then the path $i_{\ast}(\gamma_{\theta})$ rotates the plane by one whole turn, and hence also rotates the points $\{ p_1, ... , p_k \} $ by one whole turn. Therefore $i_{\ast}(\gamma_{\theta})$ is identified with $\Delta_k$ under the classical identification of $\pi_1(F_k, \ast)$ with $P_k$, where $\Delta_k$ is the generator of the centre  $Z_k$ of the pure braid group $P_k$. This implies,
$$i_{\ast}(\pi_1({\rm {Aff}}_{+}({\bf R}^2), \ast))\cong \langle i_{\ast}(g)\rangle \cong \langle \Delta_k \rangle \cong Z_k.$$
We now claim that for $k\geq 5$ the fundamental group of $E_k$ is isomorphic to the pure braid group $P_k$. 

 Consider the smooth inclusion $E_k\hookrightarrow F_k$, and note that both $E_k$ and $F_k$ are smooth submanifolds of ${\bf R}^{2k}$, because they are the complements of a collection of closed subspaces of codimension $2$.  Furthermore, $E_k$ is  a smooth submanifold of $F_k$. Similarly $F_k - E_k \hookrightarrow F_k$ is a smooth inclusion. We will show that the complement of $E_k$ in $F_k$ has dimension $k+2$. The claim $\pi_1(E_k, \ast )\cong P_k$ will follow from the general fact that when a smooth submanifold $N^n$ of dimension $n$ in a smooth manifold $M^m$ of dimension $m$ , with $m-n\geq 3$ then the inclusion $$i:(M - N)\rightarrow M$$ induces an isomorphism on the fundamental groups $$i_{\ast}:\pi_1((M - N), x)\rightarrow \pi_1(M,x)$$ so that $\pi_1((M - N), x)\cong \pi_1(M,x)$.

Notice $F_k$ is an open subset of ${\bf R}^{2k}$, which means ${\rm{dim}} (F_k)=2k$. An element of $F_k - E_k$ is a $k$-tuple of points $\{p_1, ... , p_k \}$, all of which lie on the same line in ${\bf R}^2$. First we will use the  real projective plane ${\bf R}P^2$ minus a point to parametrise free lines in ${\bf R}^2$. Let $P$ be the plane $\{ (x,y,z)\in {\bf R}^3 : z=1 \}$, for $[l]\in {\bf R}P^2$ consider the plane through the origin which is normal to $l$ and call it $l^{\perp}$. Denote the intersection $P\cap l^{\perp}$ by $l'$, then $[l]\mapsto l'$ is the homeomorphism that parametrises free lines in ${\bf R}^2$ by elements of ${\bf R}P^2$ minus the element corresponding to the $z$ axis  (and ${\rm{dim}}({\bf R}P^2)=2$).

Since the $k$-tuple $\{p_1, ... , p_k \} $ is allowed to move freely within the line $l$, it can be described as an open subset of ${\bf R}^k$. Hence ${\rm{dim}} (F_k - E_k)=k+2$, and we are left with $2k - (k+2) \geq 3 \Leftrightarrow k \geq 5$. Which implies that $\pi_1(E_k, \ast )\cong \pi_1(F_k, p)\cong P_k$ when $k\geq 5$.

Now from the exactness of the sequence of homotopy groups we obtain that $\pi_1(A_k, \ast)$ is isomorphic to $ P_k / Z_k$.
\end{proof}

\begin{figure}[htbp]
\begin{center}
\includegraphics[width=50ex]{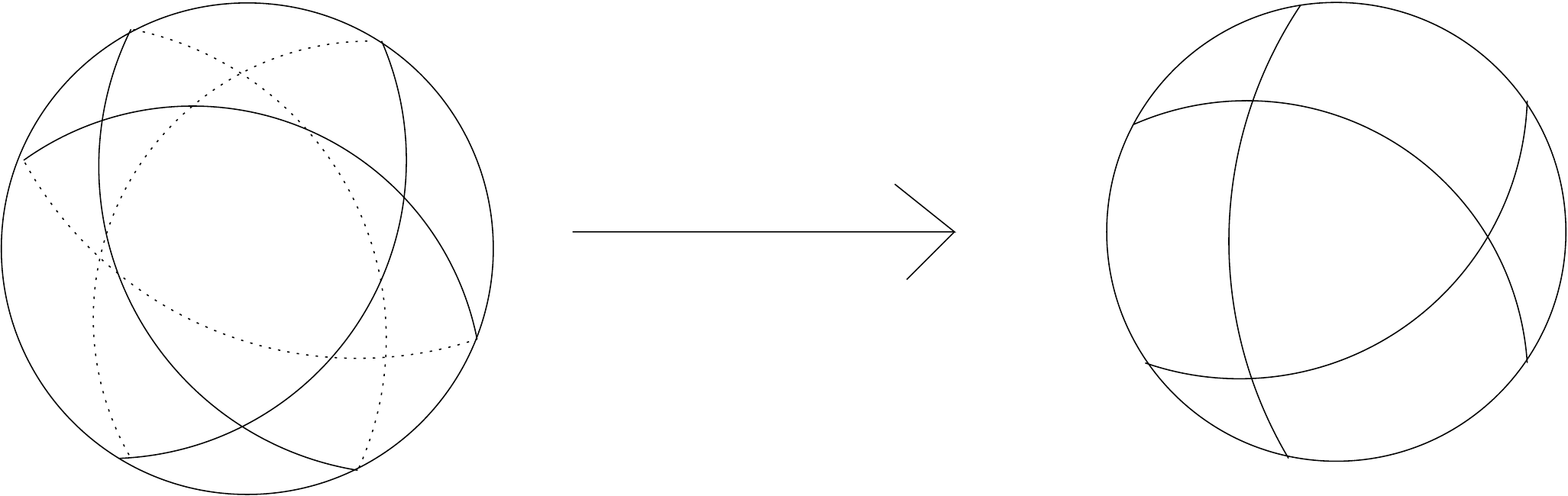}
\caption{The combinatorial structure lifts to the double orientation cover.}
\label{default}
\end{center}
\end{figure}

Notice that a presentation for $\pi_1(A_{k}, \ast)$ can be obtained from a standard presentation of $P_{k}$ adding a new relation that comes from $\Delta_{k}$.

In \cite{BMO} it was shown that $\mathbb{A}^2_{4,0}$ is homeomorphic to ${\bf R}P^2$. Taking the double orientation covering we can see $\overrightarrow{\mathbb{A}^2_{4,0}}$ is homeomorphic to $S^2$ and it has the polyhedral structure of the truncated cube.

 Then $A_4$, which in the notation of \cite{ABM} is $\overrightarrow{\mathbb{A}^2_{4,0}}(\neq)$, is homeomorphic to $S^2$ minus the vertices of the truncated cube. Because these vertices correspond to configurations where two points $p_i$ and $p_j$ coincide, i.e. $p_i=p_j$ for $i,j=1,2,3,4$. Hence $A_4$ is homeomorphic to $S^2$ minus 12 points, which correspond to the vertices of the truncated cube. Therefore the fundamental group of $A_4$ is isomorphic to the free group generated by 11 elements.

{{\bf Question: }}{\emph{ Can the cohomology of $ A_k$ be computed using the principal fibre bundle \[ {\rm{Aff}}_{+}({\bf R}^2)\stackrel{i}{\rightarrow}  E_k \rightarrow A_k \]
and spectral sequences associated to it ? }}

\end{document}